\documentclass[11pt,reqno]{article}

\usepackage[utf8]{inputenc}

\usepackage[a4paper, margin=2.5cm,top=3cm, bottom=3cm]{geometry}
\usepackage{amsmath, amsthm, amssymb, amsfonts}
\usepackage{mathtools}
\usepackage{hyperref}
\usepackage{enumitem}
\usepackage{setspace}
\usepackage{xcolor}
\usepackage{thmtools}
\usepackage{thm-restate}
\usepackage{verbatim}

\newcommand\numberthis{\addtocounter{equation}{1}\tag{\theequation}}

\theoremstyle{plain}
\newtheorem{theorem}{Theorem}[section]
\newtheorem{lemma}[theorem]{Lemma}
\newtheorem{claim}[theorem]{Claim}

\newtheorem{corollary}[theorem]{Corollary}
\newtheorem{conjecture}[theorem]{Conjecture}

\newtheorem{fact}[theorem]{Fact}
\theoremstyle{definition}

\newcommand{\Gnp}{G_{n,p}}

\newcommand{\kr}[1]{\mathrm{k}_{#1}}

\newcommand*\samethanks[1][\value{footnote}]{\footnotemark[#1]}

\title{Reconstructibility of the $K_r$-count from $n-1$ cards}

\author{
	Charlotte Knierim\thanks{Department of Computer Science, ETH Z\"urich, Switzerland\newline $\{$cknierim$\vert$anders.martinsson$\}$@inf.ethz.ch\newline
	}
	\and Anders Martinsson\samethanks[1]}

\hypersetup{
	colorlinks,
	linkcolor={red!60!black},
	citecolor={green!50!black},
	urlcolor={blue!80!black}
}

\begin{document}
\maketitle
\begin{abstract}
    The Reconstruction Conjecture of Kelly and Ulam states that any graph $G$ with $n\geq 3$ vertices can be reconstructed from the multiset $\mathcal{D}(G)$ of unlabelled subgraphs $G-v$ for all $v\in V(G)$. We refer to $\mathcal{D}(G)$ as the \emph{deck} of $G$ and $G-v\in \mathcal{D}(G)$ as the cards of $G$. This was posed in the 1940s and is still wide open today. In an effort to understand reconstructibility better, a growing collection of research is concerned with understanding what properties of $G$ can be reconstructed from a (potentially adversarially chosen) collection of $k$ cards for some $k< n$. In this paper, we show that the clique count of $G$ is reconstructible for all but one size of clique from any $n-1$ cards. We extend this result by showing that for graphs with average degree at most $3n/8-O(1)$ we can reconstruct the $K_r$-count for all $r$, and that for $r\le \log_2 n$ we can reconstruct the $K_r$-count for every graph on $n$ vertices.
\end{abstract}
\section{Introduction}
Given a graph $G$ and a vertex $v\in V(G)$, the \emph{card} $G-v$ is the unlabelled induced subgraph $G[V\setminus\{v\}]$. We call the multiset of all $n$ cards of a graph $G$ the \emph{deck} of $G$ and denote it by $\mathcal{D}(G)$. 

A natural question is whether a graph is uniquely defined by its deck, in this case we say that the graph $G$ is \emph{reconstructible}. Kelly and Ulam formulated the following conjecture in the 40s which has since become known as the \emph{Reconstruction Conjecture} \cite{kelly1942isometric,kelly1957congruence,ulam1960collection}.
\begin{conjecture}
For $n>2$, two graphs $G$ and $H$ on $n$ vertices are isomorphic if and only if $\mathcal{D}(G)= \mathcal{D}(H)$.
\end{conjecture}
This conjecture has attracted a lot of attention over the years. It is still widely open, although it has been confirmed for a few classes of graphs (e.g.\ disconnected graphs, Eulerian graphs, trees \cite{kelly1957congruence} and outerplanar graphs~\cite{giles1974reconstruction}). Moreover, it was shown by Bollob\'as~\cite{bollobas1990almost} that almost every graph can be reconstructed. This was proven by showing that the random graph $\Gnp$ is reconstructible in a very strong sense, meaning one can take as little as 3 cards (recall that there are $n$ cards in total) and still reconstruct $\Gnp$ for any $2.5\ln (n)/n \leq p \leq 1-2.5\ln (n)/n$ with probability $1-o(1)$.

As reconstructing the full graph turns out to be quite difficult, people have turned to the question of whether we can reconstruct certain \emph{properties} of a graph and how many cards are needed to do so. We say that a certain property of a given graph $G$ can be \emph{reconstructed from $k$ cards}, if, given any $k$ cards from $\mathcal{D}(G)$, it is possible to determine whether $G$ has the property or not.\footnote{In other words, a property $P$ of a given graph $G$ is reconstructible from $k$ cards if, for every graph $H$ which has at least $k$ cards in common with $G$, either both $G$ and $H$ have the property $P$, or neither of them do.} A famous lemma by Kelly \cite{kelly1957congruence} states that if we are given all the cards, then we can reconstruct the subgraph count for every proper subgraph of $G$.
It was shown that connectedness can be can be reconstructed from $\lfloor\frac{n}{2}\rfloor +2$ cards \cite{bowler2011recognizing}. This is known to be tight.

A lot of attention has been given to \emph{size reconstruction}, the reconstruction of the number of edges of the graph $G$. The first result is by Myrvold~\cite{myrvold1992degree} in 1992. She showed that the degree sequence (and thus the number of edges) can be reconstructed from $n-1$ cards for any $n\geq 7$. As further observed by Myrvold, the assumption of $n\geq 7$ is necessary, as for any $3\leq n \leq 6$, there exists pairs of $n$-vertex graphs with $n-1$ common cards, but with different edge-counts.

It took over 25 years for an improvement of this result. Brown and Fenner \cite{brown2018size} showed that $n-2$ cards suffice to determine the number of edges in $G$ for sufficiently large $n$. Recently, Groenland, Guggari and Scott~\cite{groenland2021size} showed that in fact we can miss a superconstant number of cards.

\begin{theorem}[\cite{groenland2021size}]
For $n$ sufficiently large and $k\le 0.05\sqrt{n}$, the number of edges $m$ of a graph $G$ on $n$ vertices is reconstructible from any $n-k$ cards.
\end{theorem}

They conjecture that Kelly's Lemma can be extended to hold for $n-k$ cards for any constant $k$ and $n$ large enough (see Conjecture 3.2 in \cite{groenland2021size}). Very recently, as a first step towards this Groenland et al.~\cite{groenland2021reconstructing} obtained the following result which they also use as a tool to determine the degree sequence of sparse graphs with $O(n/d^3)$ missing cards, where $d$ is the average degree of the graph.
\begin{theorem}[\cite{groenland2021reconstructing}]
\label{thrm:Groenland}
Let $d,r\in \mathbb{N}$. For any graph $G$ on $n$ vertices with average degree at most $d$, the number of cliques of size $r$ in $G$ can be reconstructed from any deck missing at most $\left(1 +\binom{2(d+1)}{r-1}\right)^{-1}(n/2-1)-d-5$ cards.
\end{theorem}
 In particular, for any fixed $r\geq 2$, their result states that the $K_r$-count can be reconstructed from $n-k$ cards for any graph with average degree $d=O(\min(n^{1/r}), (n/k)^{1/(r-1)})$. Unfortunately their result does not apply to dense graphs.
We go in a slightly different direction. Instead of asking how many cards can be missing, we explore what can be reconstructed from $n-1$ cards. In this case, Theorem~\ref{thrm:Groenland} gives the answer if $d= O(n^{1/(r-1)})$. Our first results improves this to $3n/8-O(1)$. In particular, we eliminate the dependence on the size of the clique.

\begin{restatable}{theorem}{sparse}
\label{thrm:sparse}
Let $G$ be a graph on $n\ge 7$ vertices with average degree $d\le \frac{3n}{8}-O(1)$. Then we can reconstruct the $K_r$-count of $G$ from $n-1$ cards for all $r$.
\end{restatable}

As a second result, we remove the dependence on the average degree for $r\le \log_2 n$, showing that we can count small cliques in any graph.

\begin{restatable}{theorem}{smallr}

\label{thm:smallr}
Let $G$ be a graph on $n\ge 7$ vertices and let $r\le \log_2 n$. Then we can reconstruct the $K_r$-count in $G$ from $n-1$ cards.
\end{restatable}

When proving Theorem~\ref{thrm:sparse} and Theorem~\ref{thm:smallr}, we first prove the following statement, that allows us to reconstruct the clique count for almost all $r$ for every graph. We believe that the restriction $r\ne n-\ell$ is an artefact of our proof and that it is indeed possible to reconstruct the $K_r$-count for all $r$ in any graph on $n\ge 7$ vertices.

\begin{restatable}{theorem}{main}
\label{thrm:main}
Let $G$ be a graph on $n\ge 7$ vertices and let $\ell$ be the number of maximum degree vertices in $G$. Then we can reconstruct the $K_r$-count in $G$ from $n-1$ cards for all $r\ne n-\ell$.
\end{restatable}

\paragraph{Structure of this paper}
In Section~\ref{sec:proof} we start by introducing some useful tools, then develop our main lemma, a structural result on graph classes that are `harder' in terms of reconstructing the $K_r$-count, and we then prove Theorem~\ref{thrm:main}. In Section~\ref{sec:rln} we talk about what these structural results imply for the case $n = \ell + r$. We also prove Theorems~\ref{thrm:sparse} and~\ref{thm:smallr} in this section. 
\section{Proof of Theorem \ref{thrm:main}}    
\label{sec:proof}
The $K_r$-count of a graph $G$, denoted $\kr{r}(G)$, is the number of copies of $K_r$ contained in $G$. We will frequently refer to the number of cliques containing one particular vertex. As this is extension of the classical degree notion, we define the $K_r$-degree of a vertex $v\in V(G)$, denoted $\deg_r(v,G)$, as the number of cliques of size $r$ that the vertex $v$ is contained in. In particular, $\deg_2(v,G)$ denotes the degree of $v.$ We omit the second parameter if the graph $G$ is clear from the context.

We will denote the vertices of the graph $G$ whose properties we aim to reconstruct by $v_1, \dots, v_n$ where $n=|V(G)|$. We will denote the cards visible to us by the (multi-)set $\mathcal{D}' = \mathcal{D}(G)\setminus\{G-v_h\}$, where $v_h$ is the \emph{hidden vertex} of $G.$ Given any graph $G'\subseteq G$, we will denote by $V_d(G')$ the set of vertices in $G'$ of degree $d$ in $G'$. Similarly, we will denote by $\mathcal{D}_d(G)$ and $\mathcal{D}'_d$ the set of cards in $\mathcal{D}(G)$ and $\mathcal{D}'$ where the removed vertex has degree $d$.

We start this section by giving some intuition on what are sufficient conditions to reconstruct the $K_r$-count of $G$ easily. We start with the following simple fact which follows from the observation that any non-edge is visible on at least one of any three cards.
\begin{fact}
Let $G$ be a graph on $n\ge 3$ vertices. Then we can verify whether $G$ is a clique from any $3$ cards.
\label{fact:clique}
\end{fact}
This solves the case $r=n$, where we know that $\kr{r}(G)=0$ (respectively $\kr{r}(G)=1$) if $G$ is not a clique (respectively is a clique) on $n$ vertices.

If we can identify a visible card $\mathcal{C}=G-v_i\in \mathcal{D}'$ and the corresponding $K_r$-degree of the vertex $v_i$ removed from the card, then we can easily deduce that 
\begin{equation}\label{eq:carddegrel}\kr{r}(G) = \kr{r}(\mathcal{C})+\deg_r(v_i).
\end{equation}
A natural way to approach this is to try to find another card $\mathcal{C'}\in \mathcal{D}'$ on which $v_i$ can be identified (as the vertex removed from card $\mathcal{C}$) and where $\deg_2(v_i, G)=\deg_2(v_i, \mathcal{C}')$ (i.e. the vertices removed from cards $\mathcal{C}$ and $\mathcal{C}'$ are non-neighbours in $G$). Then $\mathcal{C}'$ contains the full neighbourhood of $v_i$ and this allows us to compute $\deg_r(v_i, G)=\deg_r(v_i, \mathcal{C}')$.

The main problem is: how do we find such a pair? Given a card $\mathcal{C}=G-v_i$ and the edge-count of $G$, we can deduce $\deg_2(v_i)$ by counting how many edges are missing on $\mathcal{C}$ (all of them must be incident to $v_i$). Clearly $G$ can have a lot of vertices with equal degrees, making this identification very difficult. Additionally, if we see a vertex on a card, it is very hard to determine whether or not it is visible with its full neighbourhood or not (i.e. whether the vertex is adjacent to the vertex that was deleted to obtain the card or not). 

There is one condition to make this identification easier. Given an integer $k>0$, if there is \emph{no} vertex of degree $k+1$ in $G$, then we can be sure that, whenever we see a vertex of degree $k$ on a card, this vertex is visible with full degree. In particular, this holds for vertices of maximum degree. They will play a special role as they have the property that we can identify them when we see them on the card of a non-neighbour. We first introduce some useful tools.
\subsection{Tools}

We start this section by introducing some general tools that we repeatedly use. As we have $n-1$ cards, we repeatedly use the fact that we can reconstruct the degree sequence from these cards. 

\begin{theorem}[\cite{myrvold1992degree}]
\label{thm:myrvold}
Let $G$ be a graph on $n\ge 7$ vertices, then we can reconstruct the degree sequence of $G$ from $n-1$ cards.
\end{theorem}
In particular, as the degree sequence immediately gives us the number of edges by the handshake lemma, we can deduce, for any given card, the degree of the removed vertex  by counting the number of edges visible on the card --- the removed vertex is adjacent to all the missing edges. Note that this directly gives us $\deg_2(v_h)$. Consequently, for every $d$ we can determine from  $\mathcal{D}'$ the (multi-)set $\mathcal{D}'_d$ of visible cards whose removed vertex has degree $d$. We will not only need the degrees of the vertices but also information about the number of cliques of size $r$ a specific vertex is contained in. We see the following lemma as an extension of Theorem~\ref{thm:myrvold}.
\begin{lemma}
\label{lem:deg_seq}
Let $G$ be a graph on $n$ vertices. Suppose $\kr{r}(G)$ is known. For any $\mathcal{C}\in\mathcal{D}'$ we have
$$\deg_r(v_i)=\kr{r}(G)-\kr{r}(\mathcal{C}),$$
where $v_i$ denotes the vertex removed from card $\mathcal{C}$. Moreover, $\deg_r(v_h)$ is uniquely determined by $\mathcal{D'}$ and $\kr{r}(G)$.%In particular, this allows us to compute $\deg_r(v_h)$.
\end{lemma}
\begin{proof}
%Given $\kr{r}(G)$, for every visible card $\mathcal{C}=G-v_i$ we have $\deg_r(v_i)=\kr{r}(G)-\kr{r}(\mathcal{C})$. 
%As described above, we can get $\deg_r(v_i)$ by counting the number of $K_r$'s on the card $\mathcal{C}$.
The first statement follows from \eqref{eq:carddegrel}. Thus the only thing left to do is to calculate the $K_r$-degree of the hidden vertex. We get 
\[\kr{r}(G) = \frac{1}{r}\sum_{i\in[n]}\deg_r(v_i).\]
This implies that 
\[\deg_r(v_h)= r\cdot\kr{r}(G) - \sum_{i\in[n]\setminus \{h\}}\deg_r(v_i).\qedhere\]
\end{proof}

 We observe that, by a double counting argument, once we have figured out the $K_r$-degree of the hidden vertex, we can deduce the $K_r$-count of the graph. Note that figuring out the $K_r$-degree of the hidden vertex is not easy in general, but we will encounter cases where it easily follows from the situation we consider.
 
\begin{lemma}\label{fact:hiddencard_kr}
Let $G$ be a graph on $n\ge 7$ vertices and let $r\le n-2$. Then $\kr{r}(G)$ is uniquely determined by $\mathcal{D}'$ and $\deg_r(v_h)$.
%Then, given $deg_r(v_h)$, we can deduce $\kr{r}(G)$ from $\mathcal{C}_1,\ldots,\mathcal{C}_{n-1}$. 
\end{lemma}
\begin{proof}
Every clique of size $r$ in $G$ is present on exactly $n-r$ cards in $\mathcal{D}(G)$. Hence
\[(n-r)\kr{r}(G)= \sum_{\mathcal{C}\in\mathcal{D}(G)} \kr{r}(\mathcal{C}).\]
Substituting the clique count of the missing  card $\kr{r}(G-v_h)$ by $\kr{r}(G)-\deg_r(v_h)$ gives 
\begin{align*}
    (n-r)\kr{r}(G)&= \sum_{\mathcal{C}\in\mathcal{D}'} \kr{r}(\mathcal{C}) +(\kr{r}(G)-\deg_{r}(v_h)),
    \end{align*}
    and rearranging this equation we obtain
    \begin{align*}
    \kr{r}(G) &= \frac{\sum_{\mathcal{C}\in\mathcal{D}'} \kr{r}(\mathcal{C}) -\deg_{r}(v_h)}{n-(r+1)},
\end{align*}
where the right-hand side is completely determined by $\mathcal{D}'$ and $\deg_{r}(v_h)$ as desired.
%We can compute the last expression as we see the cards $\mathcal{C}_1,\ldots,\mathcal{C}_{n-1}$ and know $\deg_{r}(v_h)$ by assumption.
\end{proof}

The following lemma resolves the $K_r$-count if our graph has a very specific structure. Let $\Delta=\Delta(G)$ denote the maximum degree in $G$, and let $\ell=\ell(G)$ denote the number of vertices with degree $\Delta$.
\begin{lemma}
\label{lem:clique}
Let $G$ be a $n$-vertex graph with precisely $\ell$ vertices of maximum degree and $\Delta \le n-2$. If the vertices of maximum degree form a clique and the hidden vertex has maximum degree, then we can reconstruct the $K_r$-count of $G$ from $n-1$ cards unless $n-\ell-r= 0$.
\end{lemma}
\begin{proof}
 Note that, by Theorem~\ref{thm:myrvold}, the quantities $\ell$, $\Delta$, $\deg_2(v_h)$, as well as the degree of the missing vertex of any card $\mathcal{C}\in\mathcal{D}'$ are all uniquely determined by $\mathcal{D}'$. Moreover, we can deduce from $\mathcal{D}'$ whether the $\ell$ vertices of degree $\Delta$ in $G$ form a clique by checking if there is a visible card $\mathcal{C}\in\mathcal{D}'_\Delta$ that satisfies $\Delta(\mathcal{C})=\Delta.$ Thus, it is possible to determine from $\mathcal{D}'$ alone whether the conditions of the lemma are satisfied.

We will prove the lemma by using a double counting argument.
We call a pair of a clique $K$ of size $r$ and a vertex $v\in K$ an \emph{$r$-clique rooted in vertex $v$} in $G$. We count the number of rooted $r$-cliques of $G$ for which the root vertex has degree $\Delta$. (In the case where an $r$-clique contains multiple vertices with degree $\Delta$, it is counted with multiplicity.) %We count the number of all $K_r$s rooted in a vertex of degree $\Delta$ (visible with all $\Delta$ neighbours) among all cards that we see. If a $K_r$ has multiple such vertices, we count it once for each of the vertices. 
Let $\mathcal{K}$ be the set of all $r$-cliques in $G$. Clearly we have $|\mathcal{K}|=\kr{r}(G)$.

The first way we can approach the  number of rooted $r$-cliques of $G$ is to count the number of times any card $\mathcal{C}\in\mathcal{D}'$ contains an $r$-clique rooted in a vertex $v$ with $\deg_2(v, \mathcal{C})=\Delta$. This counts the $r$-cliques of $G$ rooted in a vertex $v$ with $\deg_2(v, G)=\Delta$ once for each card corresponding to a non-neighbour of $v$, or $n-1-\Delta$ times in total. Crucially, the hidden card does not correspond to a non-neighbour of $v$ as the vertices of degree $\Delta$ form a clique. This yields
\[\sum_{\mathcal{C}\in\mathcal{D}'}\sum_{w\in V_\Delta(\mathcal{C})} \deg_r(w,\mathcal{C}) =\sum_{K\in \mathcal{K}}|V_\Delta(G)\cap K|\cdot (n-1-\Delta).\numberthis\label{eq:clique} \]

On the other hand, we can obtain a closely related quantity by counting the number of (non-rooted) cliques of size $r$ on cards $\mathcal{C}=G-v_i$ for which $\deg_2(v_i,G)<\Delta$, i.e.\ $\mathcal{C}\in\mathcal{D}'\setminus\mathcal{D}'_\Delta$. Each $r$-clique of $G$ is present on $n-r$ cards in $\mathcal{D}(G)$. If we only look at cards corresponding to a vertex of degree less than $\Delta$, such a clique $K$ is now on present on $n-r-|V_\Delta(G)\setminus K|=n-r-\ell+|V_\Delta(G)\cap K|$ of these cards.
We can be sure that we see all these cards as the hidden vertex has degree $\Delta.$ Hence
%$\#\{\text{vertices of degree $\Delta$ not in this $K_r$}\}= \ell - \#\{\text{vertices in $K_r$ with degree $\Delta$}\}$.
\begin{align*}
    \sum_{\mathcal{C}\in\mathcal{D}'\setminus\mathcal{D}'_\Delta}\kr{r}(\mathcal{C}) &= \sum_{K\in \mathcal{K}}\left(n-r-\ell+|V_\Delta(G)\cap K|\right)\\
    & =  \sum_{K\in \mathcal{K}}(n-r-\ell) + \sum_{K\in \mathcal{K}}|V_\Delta(G)\cap K|\\
    & = \kr{r}(G)\cdot(n-r-\ell) + \sum_{\mathcal{C}\in\mathcal{D}'} \sum_{w\in V_\Delta(\mathcal{C})} \deg_r(w,\mathcal{C}) / (n-1-\Delta),   
\end{align*}
where the last step follows from~\eqref{eq:clique}.
Rearranging the above, we get 
\[\kr{r}(G)\cdot(n-r-\ell)= \sum_{\mathcal{C}\in\mathcal{D}'\setminus\mathcal{D}'_\Delta}\kr{r}(\mathcal{C})-\sum_{\mathcal{C}\in\mathcal{D}'} \sum_{w\in V_\Delta(\mathcal{C})} \deg_r(w,\mathcal{C}) / (n-1-\Delta),\]
 which allows us to compute the number of $K_r$'s in $G$ whenever $n-r-\ell \ne 0$. Note that in the last step we also need that $\Delta \le n-2$ to avoid dividing by zero.
\end{proof}
\subsection{Proof of the main result}
We will prove the main theorem in several steps.
First, we prove the result for $\Delta = n-1$ and $\Delta=n-2$. 

\begin{lemma}
\label{lem:n-1}
 Let $G$ be a graph on $n\ge 7$ vertices with $\Delta(G)=n-1$. Then we can reconstruct $\kr{r}(G)$ from $n-1$ cards.
\end{lemma}
\begin{proof}
First observe that for $r=n$ we just need to check whether $G$ is a clique which can be done by Fact~\ref{fact:clique} as $n-1\ge 3$. Thus in the following we assume $r\le n-1$.
Note that if $\mathcal{D}'_\Delta$ is non-empty, then we can easily reconstruct the graph $G$ by taking any such card and adding one vertex adjacent to all vertices on the card. Thus we assume that we do not have such a card, which implies that the hidden vertex is the unique vertex of maximum degree. Note that as $\Delta = n-1$ the hidden vertex is adjacent to all other vertices and thus not visible with its complete neighbourhood on any card.

Under this assumption, we determine $\kr{r}(G)$ %$\deg_r(v_h, G)$ 
by induction on $r$ for all $2\leq r \leq n-2$. For $r=2$ this is equal to the edge count and thus known. Now, suppose $r\geq 3$. By the induction assumption, we know $\kr{r-1}(G)$. Hence, by Lemma~\ref{lem:deg_seq} we also know $\deg_{r-1}(v_h, G)$ from $k_{r-1}(G)$. Given these to quantities, we claim that the $K_r$-degree of $v_h$ is given by
\[\deg_r(v_h,G)=\kr{r-1}(G)-\deg_{r-1}(v_h,G).\]
As $v_h$ is adjacent to all other vertices in the graph, the $K_r$-degree of $v_h$ is equal to the number of copies of $K_{r-1}$ in $G-v_h$, which is precisely $\kr{r-1}(G)-\deg_{r-1}(v_h, G)$. When $r\le n-2$, this allows us to use Lemma~\ref{fact:hiddencard_kr} to calculate $\kr{r}(G).$

This only leaves case when $v_h$ is the unique vertex of degree $\Delta=n-1$ and $r=n-1$. Note that if $V(G)\setminus\{v_i\}$ and $V(G)\setminus\{v_j\}$ are both cliques for $i\neq j$ it would follow that all $n-2\geq 5$ vertices in $V(G)\setminus\{v_i, v_j\}$ have degree $\Delta=n-1$. Similarly, $V(G)\setminus\{v_h\}$ cannot be a clique, as this would imply that all vertices of $G$ have degree $\Delta$. Hence, $G$ either contains no cliques of size $n-1$, or it has exactly one clique of size $n-1$, which is given by $V(G)\setminus\{v_i\}$ for some vertex $v_i\neq v_h$. This means that there is a card (which is not the hidden card) on which we can see the clique if it exists. Thus it is easy to check from the $n-1$ available cards.
\end{proof}

\begin{lemma}
\label{lem:n-2}
 Let $G$ be a graph on $n\ge 7$ vertices with $\Delta(G)=n-2$ and exactly $\ell$ vertices of maximum degree, then we can reconstruct $\kr{r}(G)$ from $n-1$ cards unless $n-r-\ell = 0$.
\end{lemma}
\begin{proof}
Note that the condition on $\Delta(G)$ implies that $\kr{n}(G)=0$, and that the only case where $\kr{n-1}(G)>0$ is if $G$ consists of a clique of size $n-1$ together with an isolated vertex, which can be deduced from the degree sequence of $G$. Thus we may assume that $r\leq n-2$.

If $v_h$ is the unique vertex with maximum degree, then we know that there is a card $\mathcal{C}\in\mathcal{D}'$ containing a vertex of degree $\Delta$ -- namely the card corresponding to the unique non-neighbour of $v_h$. But then we can easily count $\deg_r(v_h,G)$ on this card and apply Lemma~\ref{fact:hiddencard_kr} to finish the proof.

If there are multiple cards in $\mathcal{D}(G)$ corresponding to vertices of maximum degree (in particular, $\mathcal{D}'_\Delta$ is non-empty), observe the following. If any card in $\mathcal{C}\in \mathcal{D}'_\Delta$ contains a vertex of degree $\Delta$, then we can reconstruct $G$ by adding a vertex to $\mathcal{C}$ and connecting it to all vertices except for the one of degree $\Delta$. Such a card must exist whenever $G$ contains two non-adjacent vertices of degree $\Delta$. Thus the only remaining case to consider is when the vertices of maximum degree in $G$ form a clique.
If in this case the hidden vertex has degree $\Delta$, then we can use Lemma~\ref{lem:clique} to determine $\kr{r}(G)$. 

Otherwise the hidden vertex does not have degree $\Delta$, and we will proceed by induction on $r$. In this case, we prove the statement even for $n-r-\ell = 0$.

\begin{claim}
 In the case $\Delta = n-2$, we can reconstruct $\kr{r}(G)$ from $n-1$ cards if $\deg_2(v_h) < \Delta$.
\end{claim}
As in Lemma~\ref{lem:n-1} we use $r=2$ as a base case for our induction which follows immediately from Theorem~\ref{thm:myrvold}. Assume now that we know $\kr{r-1}(G)$. 
We already established in \eqref{eq:carddegrel} that for every card $\mathcal{C}$ and the corresponding removed vertex $v_i$ we have $\kr{r}(G) = \deg_r(v_i)+\kr{r}(\mathcal{C}).$ We can easily count $\kr{r}(\mathcal{C})$ for every $\mathcal{C}\in \mathcal{D}'$. For vertices of degree $\Delta$ we can also say something about $\deg_r(v_i)$. As $\Delta = n-2$ we know that a vertex $v_i$ of degree $\Delta$ has a unique non-neighbour (say $v_{j_i}$). Then we know that $v_i$ extends all $K_{r-1}$'s visible on $\mathcal{C}$ to $K_r$'s, except for the ones containing $v_{j_i}$. We get 
\[\deg_r(v_i,G)=\kr{r-1}(\mathcal{C})-\deg_{r-1}(v_{j_i},\mathcal{C}).\]

Of course we cannot hope to identify $v_{j_i}$ on any card in $\mathcal{D}'_\Delta$. But we know that given a card $G-v_{j}$ counting the number of vertices that are visible with maximum degree on this card gives us the number of degree $\Delta$ vertices that have the removed vertex of this card as the unique non-neighbour. Instead of looking at a particular card in $\mathcal{D}'_\Delta$, we thus consider the sum over all these cards, which gives
\begin{align*}
    \ell \cdot \kr{r}(G) &=
    \sum_{v_i\in V_\Delta} \deg_r(v_i) + \kr{r}(G-v_i)\\
    &= \sum_{v_i\in V_\Delta(G)} \kr{r-1}(G-v_i)-\deg_{r-1}(v_{j_i},G-v_i)+\kr{r}(G-v_i)\\
    &=\sum_{\mathcal{C}\in\mathcal{D}'_\Delta}\left(\kr{r}(\mathcal{C})+\kr{r-1}(\mathcal{C})\right) - \sum_{v_j\in V(G)} \deg_{r-1}(v_j)|V_\Delta(G)\setminus N_G(v_j)|\\
    &=\sum_{\mathcal{C}\in\mathcal{D}'_\Delta}\left(\kr{r}(\mathcal{C})+\kr{r-1}(\mathcal{C})\right) - \sum_{\mathcal{C}\in\mathcal{D}(G)}\left((\kr{r-1}(G) - \kr{r-1}(\mathcal{C}))\cdot |V_\Delta(\mathcal{C})|\right).\numberthis\label{eq:n-2} 
\end{align*}
We claim that it is possible to evaluate the right-hand side of this expression given $\mathcal{D}'$. Indeed, the contribution to this sum from each card in $\mathcal{D'}$ can be directly identified from the corresponding card. (Recall that $\kr{r-1}(G)$ is known by the induction hypothesis.) So it only remains to determine the value of
$$ (\kr{r-1}(G) - \kr{r-1}(G-v_h))\cdot |V_\Delta(G-v_h)|.  $$
The factor $\kr{r-1}(G) - \kr{r-1}(G-v_h) = \deg_r(v_h)$ is known by Lemma \ref{lem:deg_seq}. Moreover, as every vertex of maximum degree has exactly one non-neighbour, we know that it is present with its full degree on exactly one card from $\mathcal{D}(G)$. Thus we have $|V_\Delta(G-v_h)| = \ell - \sum_{\mathcal{C}\in \mathcal{D}'}|V_\Delta(\mathcal{C})|$, which again can be computed given $\mathcal{D}'.$ Putting all of this together gives us an expression for $\kr{r}(G)$ as a function of $\mathcal{D}'$, as desired.% Thus we can compute the second sum over all cards and thus evaluate $\kr{r}(G)$ from \eqref{eq:n-2}. The first sum only goes over the cards corresponding to a vertex of maximum degree, as the hidden vertex does not have maximum degree by assumption all of the needed cards are visible.
\end{proof}

 In the following we will treat the remaining case $\Delta\le n-3$. Taking all cards that belong to a maximum degree vertex, we know that there needs to be an assignment of the maximum degree vertices to the maximum degree cards (including possibly the hidden card). Assuming $\Delta\le n-3$ ensures that every vertex has at least two non-neighbours. In particular, it ensures that every vertex of maximum degree is visible on at least one card (\emph{not} the hidden card) of a non-neighbour. This allows us to identify the vertex as a vertex of maximum degree (as there are no vertices with degree $\Delta+1$). As all the $K_r$'s touching a particular vertex must be contained in its neighbourhood, this allows us to count the $K_r$-degree of vertices of maximum degree.
Let $\mathcal{A}$ be the set of all the $K_r$-degrees we see on vertices of maximum degree. More formally,
\[\mathcal{A} := \left\{\deg_r(w,G)| w\in V_\Delta(G)\right\} = \left\{\deg_r(w,\mathcal{C})| \mathcal{C}\in\mathcal{D}', w\in V_\Delta(\mathcal{C})\right\}.\]

If the hidden vertex does not have maximum degree, then we see all the cards corresponding to maximum degree vertices. As we know that $\kr{r}(G) = \deg_r(v_i)+\kr{r}(G-v_i)$ for all these cards, we know that the card $\mathcal{C}\in\mathcal{D}'_\Delta$ with the lowest $k_r(\mathcal{C})$ needs to go together with a vertex in $V_\Delta(G)$ of the highest $K_r$-degree and the card with the highest $k_r(\mathcal{C})$ needs to go with a vertex of the lowest $K_r$-degree. Then we can determine the assignment of vertices of maximum degree to the corresponding cards up to permutation of vertices with the same $K_r$-degree. Even if the hidden card belongs to a vertex of maximum degree, we know that at least one of these pairings must happen. The following lemma formalises this. 

\begin{lemma}
\label{lem:two_choices}
Let $G$ be a graph on $n\ge 7$ vertices with $\Delta(G)\le n-3$ and let $r\ge 2 $ be an integer. Given any $n-1$ cards such that $v_h$ is not the unique vertex of maximum degree in $G$, at least one of the following equations is always true 
\begin{align*}\kr{r}(G) = 
\min \mathcal{A} +& \max_{\mathcal{C}\in\mathcal{D}'_\Delta}\kr{r}(\mathcal{C}),\numberthis \label{eq:minmax}\\&\text{or}\\
\kr{r}(G)=\max \mathcal{A} +& \min_{\mathcal{C}\in\mathcal{D}'_\Delta}\kr{r}(\mathcal{C}).\numberthis\label{eq:maxmin}
\end{align*}
 If $\deg_2(v_h)<\Delta$ or the $K_r$-degree of the hidden vertex is not unique among $V_\Delta(G)$, then \eqref{eq:minmax} and \eqref{eq:maxmin} are both true and we can precisely determine $\kr{r}(G)$.
\end{lemma}

\begin{proof}[Proof of Lemma~\ref{lem:two_choices}]
Given $n-1$ cards and assuming $\Delta\leq n-3$, we can determine the set $\mathcal{A}$ of $K_r$-degrees of the vertices in $G$ with degree $\Delta$. We also have all but at most one of the cards where the removed vertex has maximum degree in $\Delta$. Let
$$\mathcal{B}:=\{\kr{r}(\mathcal{C})|\mathcal{C}\in\mathcal{D}'_\Delta\},$$
and observe that
\[\mathcal{A}\supseteq \kr{r}(G)-\mathcal{B},\numberthis\label{eq:AB}\]
where $\kr{r}(G)-\mathcal{B}$ denotes the set of all integers $\kr{r}(G)-b$ for all $b\in\mathcal{B}$, and where the right-hand side equals $\mathcal{A}\setminus\{\deg_r(v_h)\}$ if $\deg_2(v_h)=\Delta$ and the $K_r$-degree of $v_h$ is unique among all vertices in $V_\Delta(G)$, and equals $\mathcal{A}$ otherwise. Note that, by assumption, $\mathcal{B}$ is non-empty. Since the left-hand and right-hand sides of~\eqref{eq:AB} differ by at most one element, they must either have the same minimiser, which implies~\eqref{eq:minmax}, or the same maximiser, which implies~\eqref{eq:maxmin}.

If we have equality in~\eqref{eq:AB} then they have the same minimiser and maximiser and both~\eqref{eq:minmax} and~\eqref{eq:maxmin} are true. In this case they have to give the same value so we have precisely determined $\kr{r}(G)$.
\end{proof}

For convenience of the reader, for the rest of the proof we will always assume that $\ell$ is the number of vertices of maximum degree in $G$, and that the vertices of $G$, $v_1, \dots, v_n$ are ordered lexicographically decreasing by the values of $(\deg_2(v_i), \deg_r(v_i))$, with ties broken arbitrarily. In particular, $v_1, \dots, v_\ell$ denotes the vertices of maximum degree with $$\deg_r(v_1)\geq \deg_r(v_2) \geq \dots \geq \deg_r(v_\ell),$$
where $\deg_{r}(v_1)=\max \mathcal{A}$ and $\deg_{r}(v_\ell)=\min \mathcal{A}$. Analogously, we will order the visible cards $$\mathcal{C}_1, \dots, \mathcal{C}_{n-1}\in \mathcal{D}'$$
lexicographically increasing with respect to $(e(\mathcal{C}_i), \kr{r}(\mathcal{C}_i))$. In other words, depending on whether or not $v_h$ has degree $\Delta$, the visible cards corresponding to maximum degree vertices are either $\mathcal{C}_1, \mathcal{C}_2, \dots, \mathcal{C}_{\ell-1},$ or $\mathcal{C}_1, \mathcal{C}_2, \dots, \mathcal{C}_{\ell}$, ordered decreasing by their $K_r$-count, i.e. $$\kr{r}(\mathcal{C}_1)\le \kr{r}(\mathcal{C}_2)\le \ldots \le\kr{r}(\mathcal{C}_{\ell-1}) (\le \kr{r}(\mathcal{C}_\ell)).$$

\begin{lemma}
\label{lem:main}
Let $G$ be graph on $n\ge 7$ vertices with maximum degree $\Delta$ and let $r\ge 2$ be a constant. Then, given any $n-1$ cards, either
\begin{enumerate}
    \item we can determine $\kr{r}(G)$, or
    %\item all $K_r$-degrees of the vertices $w_1,\ldots,w_\ell$ are distinct.
    \item we have that the $K_r$-degrees of the vertices of maximum degree are all distinct, $\ell\geq 2$, the hidden vertex has degree $\Delta$, and one of the two assignments
    $\mathcal{C}_i=G-v_i\; \forall i\in [\ell-1],$
    or
    $\mathcal{C}_i=G-v_{i+1}\; \forall i\in [\ell-1]$
    must hold.
\end{enumerate}
\end{lemma}
\begin{proof} 
For now, assume $\Delta \le n-3$. We will argue about why the statement is true in the cases $\Delta = n-1$ and $\Delta = n-2$ at the end of this proof.

In this proof, we will assume that the hidden vertex $v_h$ has degree $\Delta$ and that the $K_r$-degree of this vertex is unique. If either of these was not the case, we could use Lemma~\ref{lem:two_choices} to determine $\kr{r}(G)$. If the hidden vertex is the unique vertex of maximum degree, then we can see it on at least one card and thus know its $K_r$-degree. By Lemma~\ref{fact:hiddencard_kr} we can determine $\kr{r}(G)$.
Note that the degrees of the maximum degree vertices are distinct if and only if $|\mathcal{A}|=\ell$. We make a case distinction on whether this is the case or not.
\paragraph{\textbf{Case I: $|\mathcal{A}|=\ell$}}
In this case, all the $K_r$-degrees of the vertices of maximum degree are distinct. This means that $\deg_r(v_1)>\deg_r(v_2)>\ldots>\deg_r(v_\ell)$. Lemma~\ref{lem:two_choices} tells us that there are two possible values for $\kr{r}(G)$: $ \deg_r(v_1)+\kr{r}(\mathcal{C}_1)$ and $\deg_r(v_\ell)+\kr{r}(\mathcal{C}_{\ell-1})$. As any pair of a card $\mathcal{C}_i$, $i\in [\ell-1]$ and the corresponding removed vertex $v_j$, $j\in[\ell]$ satisfies $\kr{r}(\mathcal{C}_i)+\deg_r(v_j, G)=\kr{r}(G)$, this gives us two potential ways to map cards to vertices. Think of this as a bipartite graph with the vertices of maximum degree on one side and the visible cards corresponding to vertices of maximum degree on the other side. A vertex $v_i$ is connected to a card $\mathcal{C}_j$ with a red edge if $\deg_r(v_i)+\kr{r}(\mathcal{C}_j) = \deg_r(v_1)+\kr{r}(\mathcal{C}_1)$ and with a blue edge if $\deg_r(v_i)+\kr{r}(\mathcal{C}_j) = \deg_r(v_\ell)+\kr{r}(\mathcal{C}_{\ell-1})$. 

If either the red or the blue edges do not assign each card $\mathcal{C}_1,\dots, \mathcal{C}_{\ell-1}$ to a unique vertex, then we can conclude that the corresponding candidate value for $\kr{r}(G)$ is invalid and we can determine $\kr{r}(G)$ as the remaining option.

This means we have determined $\kr{r}(G)$ in all cases except those where both colours lead to a valid matching. If both colours lead to a valid matching, then we have that the true assignment is either $\forall i\in [\ell-1] \colon \mathcal{C}_i = G-v_i$ (using the red edges) or $\forall i\in [\ell-1] \colon \mathcal{C}_i =G-v_{i+1}$ (using the blue edges). A visualisation of this can be seen in Figure~\ref{fig:2_options}.

\begin{figure}
    \centering
    \includegraphics[scale=0.7]{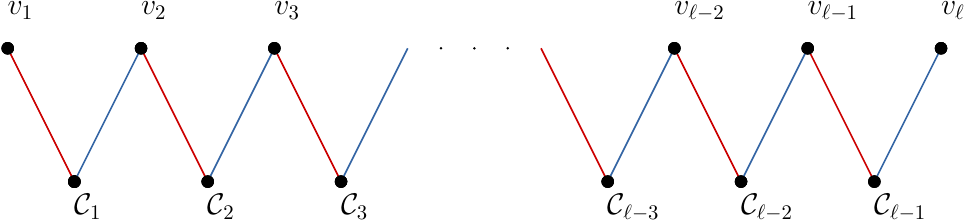}
    \caption{Illustration of Case $(2)$ of Lemma \ref{lem:main}. Given $n-1$ cards that do not uniquely determine $\kr{r}(G)$, there are only two possibilities for which vertex of $G$ could be the removed vertex from the cards in $\mathcal{D}'_\Delta$.}\label{fig:2_options}
\end{figure}

\paragraph{\textbf{Case II: $|\mathcal{A}|<\ell$}}
In this case, there are two vertices of maximum degree with the same $K_r$-degree. We divide this case into two subcases. In both subcases we want to find the largest $a\in \mathcal{A}$ that appears as the $K_r$-degree of multiple vertices in $V_\Delta(G)$. Then we show that, given this $a$, we can determine $\kr{r}(G)$. %At first, we will look at the case $\Delta= n-3$.

At first, we look at the case $\Delta=n-3$. Note that each of the maximum degree vertices is visible with their full degree on $1$ or $2$ cards in $\mathcal{D}'$ (depending on whether it is adjacent to the hidden vertex). Moreover, as $\deg_2(v_h)=\Delta=n-3$, there are at most two vertices in $V_\Delta(G)$ that are visible with their full degree on only one card in $\mathcal{D}'$ (the non-neighbours of $v_h$). We can determine the actual number of such vertices by computing $\sum_{\mathcal{C}\in\mathcal{D}'}|V_\Delta(\mathcal{C})|$. If this sum is $2\ell$ or $2\ell-1$, then we know that all but at most one vertex of $V_\Delta(G)$ is visible with its full degree on $2$ cards in $\mathcal{D}'$. In this case, a value $a$ is the $K_r$-degree of more than one vertex in $V_\Delta(G)$ if and only if among the visible cards $\mathcal{D}'$ we see at least three vertices with degree $\Delta$ and $K_r$-degree $a$. As $|\mathcal{A}|<\ell$ this happens for at least one value in $\mathcal{A} $ and we can easily identify the largest value for which this happens. 

If the aforementioned sum is equal to $2\ell-2$, then we know that both non-neighbours of $v_h$ have degree $\Delta$. If we have a visible card $\mathcal{C}_i$, for some $1\leq i \leq \ell-1$ (a card corresponding to a vertex of maximum degree), that contains two vertices of degree $\Delta$, then we can reconstruct $G$ from this card by adding a vertex and connecting it to all but the degree $\Delta$ vertices. Hence, every vertex in $V_\Delta(G)\setminus \{v_h\}$ has at most one non-neighbour in $V_\Delta(G)$, which means they have at least one (one or two) non-neighbours in $V(G)\setminus V_\Delta(G).$ In other words, $v_h$ is the unique degree $\Delta$ vertex which is not visible with its full degree on any of the cards $\mathcal{C}_{\ell},\dots, \mathcal{C}_{n-1}$, and consequently we can determine the $K_r$-degree of $v_h$ as
$$\{\deg_r(v_h)\}=\mathcal{A}\setminus \bigcup_{i=\ell}^{n-1} \{ \deg_r(w, \mathcal{C}_i) | w\in V_\Delta(\mathcal{C}_i)\}.$$
Then we have the $K_r$-degree of the hidden vertex, which by Lemma~\ref{fact:hiddencard_kr} means that we can determine $\kr{r}(G)$.

Assume now that $\Delta\le n-4$. Every maximum degree vertex is visible with its full degree on $n-\Delta-1$ or $n-\Delta-2$ cards. We know that a given $K_r$-degree is not unique if and only if the value appears more than $n-\Delta-1$ times on the visible cards as $2\cdot (n-\Delta-2)> n-\Delta -1$. 
We can easily determine the largest $a\in \mathcal{A}$ for which this is the case.

Now we want to argue that, given the largest $a\in \mathcal{A}$ such that there are multiple vertices of maximum degree with $K_r$-degree $a$, we can determine $\kr{r}(G)$. Recall that we already know by Lemma \ref{lem:two_choices} how to determine $\kr{r}(G)$ if the $K_r$-degree of $v_h$ is not unique, so the only case left to consider is if the $K_r$-degree of $v_h$ is unique. But this means that all cards corresponding to vertices with degree $\Delta$ and $K_r$-degree $a$ are visible, and we can determine the corresponding $K_r$-count as the smallest value $b$ of $\kr{r}(C_i)$ that is attained for multiple $i\in[\ell-1].$ Consequently, we have $\kr{r}(G)=a+b$, as desired.

%Now we want to argue what happens
It only remains to consider the in the cases where $\Delta = n-1$ or $\Delta = n-2$.
Recall that the case of $\Delta=n-1$ was already fully resolved in Lemma~\ref{lem:n-1} and we can always determine $\kr{r}(G)$ in this case. Let us have a look at a case when our results above cannot determine $\kr{r}(G)$ for $\Delta = n-2$. Following the proof of Lemma~\ref{lem:n-2}, we see that $\ell\geq 2$, the missing vertex has degree $\Delta$, and no card in $\mathcal{D}'_\Delta$ contains a vertex of degree $\Delta$ meaning $V_\Delta(G)$ forms a clique in $G$. In particular, this means that every vertex in $V_\Delta(G)$ is visible with its full degree on exactly one of the cards $\mathcal{D}'\setminus\mathcal{D}'_\Delta$ (the card corresponding to the unique non-neighbour of the degree $\Delta$ vertex). This lets us determine the $K_r$-degree for all the vertices in $V_\Delta(G)$. In particular, we can determine whether they are all distinct or when this is not the case easily determine the largest value that appears twice. Thus continue as in Case I or Case II respectively. This gives us that the conclusion of Lemma~\ref{lem:main} applies also for $\Delta=n-1$ or $n-2$.
\end{proof}

We want to conclude that we can determine $\kr{r}(G)$ if the maximum degree vertices do not form a clique. 

\begin{lemma}
Let $G$ be a graph on $n\ge 7$ vertices with maximum degree $\Delta\leq n-3$. We can determine $\kr{r}(G)$ if one of the cards $\mathcal{C}\in\mathcal{D}'_\Delta$ contains a vertex of degree $\Delta$.
\label{claim:clique}
\end{lemma}

\begin{proof}
\begin{comment}
\begin{figure}[t]
    \centering
    \includegraphics[scale=0.8]{figures/clique_step1.pdf}
    \caption{If we see the vertex $w_i$ on card $\mathcal{C}_j$ then we need to see $w_j$ on $\mathcal{C}_i$ and $w_{j+1}$ on $\mathcal{C}_{i-1}$ or we can determine $\kr{r}(G)$.}
    \label{fig:clique_s1}
\end{figure}
\end{comment}
Let us assume, towards a contradiction, that there exists another graph $G'$ on the same vertex set with vertices ordered analogously, and such that $\mathcal{D}'\subset \mathcal{D}(G')$  but where the graphs have different $K_r$ counts, say, $\kr{r}(G)>\kr{r}(G')$. Applying Lemma~\ref{lem:main} and observing that $\mathcal{D}'$ does not uniquely determine the $K_r$-count gives us that both $G$ and $G'$ must satisfy case $(2)$ of the statement. In particular, as $\mathcal{A}$ is uniquely determined by $\mathcal{D}'$, we must have
$$\deg_r(v_1, G)=\deg_r(v_1, G') > \deg_r(v_2, G)=\deg_r(v_2, G') > \dots > \deg_r(v_\ell, G)=\deg_r(v_\ell, G').$$
Moreover, as a pairing of a vertex and a card determines the $K_r$-count of the graph, and $\kr{r}(G)>\kr{r}(G')$, we conclude that
$$\mathcal{C}_i = G-v_i = G'-v_{i+1}\quad\forall\ell\in[\ell-1].$$

Let $H$ be the bipartite graph on vertices $a_1, \dots, a_{\ell}, b_1, \dots, b_{\ell-1}$ where $a_i$ is connected to $b_j$ by an edge if and only if $v_i$ is visible with degree $\Delta$ on $\mathcal{C}_j$, that is,  if $\mathcal{C}_j$ contains a vertex $w$ with $\deg_2(w, \mathcal{C}_j)=\Delta$ and $\deg_r(w, \mathcal{C}_j)=\deg_r(v_i)$. To see that these two formulations are equivalent, note that as $w$ has degree $\Delta$ on the card, its $K_r$-degree on the card is identical to its $K_r$-degree in $G$, and the only vertex in $G$ with this combination of degree and $K_r$-degree is $v_i$. By the assumption of the lemma, $H$ has at least one edge.

Let us now consider what conclusions we can draw from the fact that the cards $\mathcal{C}_1, \dots \mathcal{C}_{\ell-1}$ can be generated from graphs $G$ and $G'$ as above. First, from the point of view of $G$, we have that $v_i$ is visible with degree $\Delta$ on card $\mathcal{C}_j=G-v_j$ if and only if $i\neq j$ and $v_i$ and $v_j$ are not adjacent in $G$. As adjacency is symmetric, it follows that $v_j$ is visible with degree $\Delta$ on $\mathcal{C}_i$. In other words, 
\begin{equation}\label{eq:H1}a_i \not\sim b_i\quad\forall i\in[\ell-1],\end{equation}
and
\begin{equation}\label{eq:H2} a_i \sim b_j \iff a_j \sim b_i \quad \forall i,j\in[\ell-1].\end{equation}

Analogously, from the point of view of $G'$, we have that $v_{i+1}$ is visible with degree $\Delta$ on $\mathcal{C}_j=G'-v_{j+1}$ if and only if $i+1\neq j+1$ and $v_{i+1}$ and $v_{j+1}$ are not adjacent in $G'$. By symmetry, this also implies that $v_{j+1}$ is visible with degree $\Delta$ on $\mathcal{C}_{i}$. In other words, \begin{equation}\label{eq:H3}a_{i+1} \not\sim b_i\quad\forall i\in[\ell-1],\end{equation}
and
\begin{equation}\label{eq:H4} a_{i+1} \sim b_j \iff a_{j+1} \sim b_i \quad \forall i,j\in[\ell-1].\end{equation}

Consider an edge $a_i \sim b_j$ of $H$ that minimises $|i-j|$. Without loss of generality, we may assume that $i\geq j$, because if $i<j$, then $i, j\in[\ell-1]$ and by~\eqref{eq:H2} we get that $a_j \sim b_i$ is also an edge. 

If $i-j\geq 2$, then applying~\eqref{eq:H4} with $i'=i-1$ and $j'=j$ gives us that $a_{j+1}\sim b_{i-1}$ where $|(j+1)-(i-1)|=|i-j-2|<|i-j|$, which cannot be as $|i-j|$ is minimal by assumption. Similarly, we cannot have $i-j=0$ or $1$ as this contradicts~\eqref{eq:H1} or~\eqref{eq:H3} respectively. But this leaves no remaining options for the value of $i-j$, a contradiction. Hence our original assumption that there exists such a $G'$ must be false. Thus $\mathcal{D}'$ uniquely determines $\kr{r}(G)$.
\end{proof}

All that is left to do is to assemble the pieces to give a proof of Theorem~\ref{thrm:main}. We restate the theorem for convenience of the reader.

\main*

\begin{proof} Lemmas \ref{lem:n-1} and \ref{lem:n-2} give us the statement when $\Delta=n-1$ or $n-2$.
%Applying Lemma~\ref{lem:n-1} (and in the case $r=n$ we use Fact~\ref{fact:clique}) and Lemma~\ref{lem:n-2}
This leaves us with the case that $\Delta\le n-3$. By Lemma \ref{lem:two_choices}, $\kr{r}(G)$ is reconstructible from $n-1$ cards unless $v_h$ has degree $\Delta$. By Lemma \ref{claim:clique}, $\kr{r}(G)$ is reconstructible from $n-1$ cards unless the vertices of degree $\Delta$ form a clique. If both of these exceptions hold, then Lemma \ref{lem:clique} lets us compute $\kr{r}(G)$ in all cases except where $n-\ell-r=0$.
%Using Lemma~\ref{lem:two_choices}, Lemma~\ref{lem:main} and Lemma~\ref{claim:clique} we get that we can determine the $K_r$-count in almost all cases. The only case that is left, is when $v_h$ is of maximum degree and the vertices of maximum degree form a clique. In this case we use Lemma~\ref{lem:clique} to finish the proof.
\end{proof}

\section{The case $r+\ell = n$}
\label{sec:rln}

In this section we discuss the behaviour of the case $r+\ell = n$. Note that in all cases where this restriction occurs it stems from Lemma~\ref{lem:clique}. Our aim is to argue that whenever this machinery fails to uniquely determine $\kr{r}(G)$ for some fixed $r$, then $G$ must be of a special form. The assumptions from Theorems~\ref{thrm:sparse} and~\ref{thm:smallr} will allow us to uniquely determine $\kr{r}(G)$. For convenience of notation, we will again assume that the vertices of $G$ and cards in $\mathcal{D}'$ are ordered as in the previous section. That is, we assume that the vertices of $G$, $v_1, \dots, v_n$, are ordered decreasingly by $\deg_2(v_i, G)$, with vertices with equal degree ordered decreasingly by $\deg_r(v_i, G)$, and the visible cards $\mathcal{C}_1, \dots, \mathcal{C}_{n-1}$ are ordered increasingly by $e(\mathcal{C}_i)=\kr{2}(\mathcal{C}_i)$, with cards with the same edge count ordered increasingly by $\kr{r}(\mathcal{C}_i)$.

We start by observing that the maximum degree vertices must follow a distinct pattern or we can determine $\kr{r}(G)$.

\begin{lemma}
\label{lem:AP}
Let $G$ be a graph on $n\ge 7$ vertices. We can either determine $\kr{r}(G)$ from $\mathcal{D}'$, or
\begin{enumerate}
    \item $G$ contains $\ell\geq 2$ vertices of degree $\Delta\leq n-2$,
    \item $r=n-\ell$,
    \item $\deg_r(v_1) > \deg_r(v_2) > \dots > \deg_r(v_\ell)$ form an arithmetic progression, and
    \item  $\deg_{r'}(v_i) = \deg_{r'}(v_j)$ for all $r'\ne r$ and all $i, j\in [\ell]$.
\end{enumerate}
\end{lemma}
\begin{proof}
Assume $\kr{r}(G)$ is not uniquely determined by $\mathcal{D}'$. By Lemma \ref{lem:n-1}, we get that $\Delta\leq n-2$. By Theorem \ref{thrm:main}, we get $r=n-\ell$, and that $\kr{r'}(G)$ is uniquely determined by $\mathcal{D'}$ for all $r'\neq r$. To get the remaining conditions, we consider Lemma \ref{lem:main}. This implies that the $\ell\geq 2$ vertices of maximum degree $v_1,\ldots,v_\ell$ satisfy 
$$\deg_r(v_1) > \deg_r(v_2) > \dots > \deg_r(v_\ell)$$
and either
$$\mathcal{C}_i=G-v_i\; \forall i\in [\ell-1],$$
    or
    $$\mathcal{C}_i=G-v_{i+1}\; \forall i\in [\ell-1].$$ We can think of these two assignments as the red and the blue assignment as depicted in Figure~\ref{fig:2_options}.

Observe that, as $\Delta\leq n-2$, each of the vertices $v_1, \dots, v_\ell$ are visible with full degree on at least one card. As these vertices can be uniquely identified by the fact that they have degree $\Delta$ and have distinct $K_r$-degrees, it follows that we can determine $\deg_{r'}(v_i)$ for all $i\in[\ell]$ and all $r'\geq 2$ from $\mathcal{D'}.$

Given this, the two possible values of $\kr{r}(G)$ are the value given by the red assignment
$$ \kr{r}(\mathcal{C}_1)+\deg_r(v_1)=\kr{r}(\mathcal{C}_2)+\deg_r(v_2)=\dots=\kr{r}(\mathcal{C}_{\ell-1})+\deg_r(v_{\ell-1}),$$
and the value given by the blue assignment
$$ \kr{r}(\mathcal{C}_1)+\deg_r(v_2)=\kr{r}(\mathcal{C}_2)+\deg_r(v_3)=\dots=\kr{r}(\mathcal{C}_{\ell-1})+\deg_r(v_{\ell}).$$
Where, by assumption that $\kr{r}(G)$ is not uniquely determined by $\mathcal{D}'$, both expressions need to be valid. In particular, taking the difference between the equations yields
$$\deg_r(v_2)-\deg_r(v_1)=\deg_r(v_3)-\deg_r(v_2)=\dots = \deg_r(v_\ell)-\deg_r(v_{\ell-1}).$$
In other words, the $K_r$-degrees of $v_1, \dots, v_\ell$ form an arithmetic progression.

By the same line of reasoning, we can use the red and blue assignment to compute $\kr{r'}(G)$ for $r'\neq r$. However, in this case, $\kr{r'}(G)$ \emph{is} determined by $\mathcal{D}'$, meaning that values given by the red assignment $\kr{r'}(\mathcal{C}_i)+\deg_{r'}(v_i)$ and $\kr{r'}(\mathcal{C}_i)+\deg_{r'}(v_{i+1})$ need to be identical, which implies $\deg_{r'}(v_i)=\deg_{r'}(v_{i+1})$ for all $i\in[\ell-1]$, which concludes the proof of the lemma.
\end{proof}

Similar to Lemma~\ref{lem:two_choices} we can make the same argument if we can identify vertices of some given degree on cards. This happens if we have a degree $a$ such that there is no vertex of degree $a+1$. 

\begin{lemma}
Let $G$ be a graph on $n\ge 7$ vertices. If the degree sequence of $G$ has a hole, that is if there is a value $a$ such that $\delta(G)< a< \Delta(G)$ and there is no vertex of degree $a$ in $G$, then we can determine $\kr{r}(G)$ from $n-1$ cards for all $r$.
\label{lem:hole}

\end{lemma}
\begin{proof}
Assume that there is a value $a$ such that $\delta(G)< a< \Delta(G)$ and there is no vertex of degree $a$. If there are multiple such values, choose the smallest one. Note that this implies that there is a vertex of degree $a-1$. As $a-1 < n-2$ we can see every vertex with degree $a-1$ on at least one card. Clearly, when we see a vertex with degree $a-1$  on a card $\mathcal{C}\in\mathcal{D}'$ we know that this vertex is visible with all its neighbours (because if it was adjacent to the removed vertex $w$ then it would have been of degree $a$ with is not possible by assumption). By the same logic as for the vertices of maximum degree, we can thus determine the set $\mathcal{A}' := \{\deg_r(w): w\in V_{a-1}(G)\}$ from $\mathcal{D}'$.

As we know that the hidden vertex has degree $\Delta$ it cannot have degree $a-1$. By the same argument as in Lemma~\ref{lem:two_choices} we get.
Then we have \[\kr{r}(G) = \min \mathcal{A}' + \max_{\mathcal{C}\in\mathcal{D}'_{a-1}} \kr{r}(\mathcal{C}).\]

\end{proof}

Let $\omega(G)$ denote the size of the largest clique in $G$.

\begin{corollary}

Let $G$ be a graph on $n\ge 7$ vertices. If $\omega(G) < \frac{n}{2}$, then we can determine $\kr{r}(G)$ from $n-1$ cards for all $r$.
\label{fact:Deltan2}
\end{corollary}
\begin{proof} First note that it is possible to see given any $n-1$ cards whether $\omega(G)<n/2$, as any clique of size $\lceil n/2\rceil$, it is visible on all but $\lceil n/2\rceil < n-1$ cards.

Assume $\kr{r}(G)$ is not uniquely determined by $\mathcal{D}'$ for some graph $G$ with $\omega(G)<n/2.$ By Lemma \ref{lem:AP} it follows that $G$ contains a clique on $n-r$ vertices. (Namely the set of $\ell=n-r$ vertices with degree $\Delta$.) This implies that $n-r < n/2$, or, equivalently, $r>n/2$. However, then it cannot be that $\kr{r}(G)$ is not uniquely determined by $\mathcal{D}'$, as $\omega(G)<n/2$ implies that $\kr{r}(G)=0.$
\end{proof}

We can now prove Theorem~\ref{thm:smallr} and Theorem~\ref{thrm:sparse}. For the convenience of the reader we restate them here.
\smallr*
\begin{proof} Let $G$ be a graph on $n\geq 7$ vertices such that $\kr{r}(G)$ is not uniquely determined by $n-1$ cards. By Lemma \ref{lem:AP}, $n=\ell+r$, the vertices $v_1, \dots, v_\ell$ form a clique and
$$ \deg{r}(v_1)>\deg{r}(v_2)>\dots > \deg{r}(v_\ell). $$
Observe that this means that no two vertices $v_i, v_j$ for $1\leq i < j \leq \ell$ have the same neighbors among $\{v_{\ell+1}, \dots, v_n\}$. But there are only $\binom{r}{\Delta - \ell +1}$ %\le \binom{r}{\floor*{r/2}}$
choices for such a neighbourhood. % for a vertex in $\mathcal{L}$. 
 By the pigeonhole principle, it follows that $\ell=n-r\leq \binom{r}{\Delta - \ell +1}$, or, equivalently, $\binom{r}{\Delta - \ell +1}+r\geq n$. By the binomial theorem, $\binom{r}{\Delta - \ell +1}+r <\sum_{i=0}^r \binom{r}{i} \leq 2^r$, which implies that $n< 2^r$. Hence, this cannot happen if $r\leq \log_2 n.$
\end{proof}

\sparse*
\begin{proof} Let $G$ be a graph on $n\geq 7$ vertices such that $\kr{r}(G)$ is not determined by $n-1$ cards. By Corollary \ref{fact:Deltan2} we have $\omega(G)\geq n/2$, and by Lemma \ref{lem:hole} there are no holes in the degree sequence.

As $G$ has a clique of size at least $n/2$, we have that $$\deg_2(v_i)\geq \lceil n/2\rceil -1$$ for all $1\leq i\leq \lceil n/2\rceil$, and as there are no holes in the degree sequence of $G$ it follows that
$$ \deg_2(v_{i}) \geq 2\lceil n/2\rceil -1 -i$$
for all $\lceil n/2\rceil + 1\leq i \leq n.$ This implies that
\begin{align*}
    \sum_{i=1}^n\deg_2(v_i) &\geq \lceil n/2\rceil\cdot(\lceil n/2\rceil-1) + (n-\lceil n/2\rceil)\cdot \frac{(\lceil n/2\rceil -2)+(2\lceil n/2\rceil -1 -n)}2\\
    &=\frac{n}2\cdot\frac{n}2 + \frac{n}2\cdot \frac{n}{4}-O(n)=\frac{3n^2}8-O(n),
\end{align*}
which implies that the average degree $d$ of $G$ is at least $\frac{3n}8-O(1).$
\end{proof}

\section*{Acknowledgements}
We thank Ulysse Schaller and Pascal Su for helpful discussions. We also thank the anonymous referees for their careful reading of the paper and many insightful and constructive comments.

\bibliographystyle{abbrv}
\bibliography{sources}
\end{document}